\documentclass[12pt]{article}
\usepackage{amssymb,amsmath,amsfonts,amsthm,amscd,latexsym,verbatim,graphics,epsfig,indentfirst,shuffle}
\usepackage{geometry}
\def\Span{\mathrm{Span}}
\def\RB{\mathrm{RB}}
\def\Var{\mathrm{Var}}
\def\As{\mathrm{As}}
\def\Lie{\mathrm{Lie}}

\begin{document}

\begin{center}
{\Large
Embedding of post-Lie algebras into postassociative algebras}

V. Gubarev
\end{center}

\begin{center}
In honour of 80th anniversary of Professor Leonid Arkad'evich Bokut' 
\end{center}

\begin{abstract}
Applying Gr\"{o}bner---Shirshov technique, we prove that any post-Lie algebra 
injectively embeds into its universal enveloping postassociative algebra.

\end{abstract}

\section*{Introduction}

Linear operator $R$ defined on an algebra $A$ over the field $\Bbbk$
is called a Rota---Baxter operator (RB-operator, for short) of a weight $\lambda\in\Bbbk$
if it satisfies the relation
\begin{equation}\label{RB}
R(x)R(y) = R( R(x)y + xR(y) + \lambda xy), \quad x,y\in A.
\end{equation}
In this case, an algebra $A$ is called Rota---Baxter algebra (RB-algebra).

G.~Baxter defined the notion of what is now called Rota---Baxter operator on 
a (commutative) algebra in 1960 \cite{Baxter60},
solving an analytic problem. The relation~\eqref{RB} with $\lambda = 0$
appeared as a generalization of integration by parts formula.
G.-C.~Rota~\cite{Rota68}, P.~Cartier \cite{Cartier72} and others studied 
different combinatorial properties of RB-opera\-tors and RB-algebras.
In 1980s, the deep connection between Lie RB-algebras and Yang---Baxter equation
was found \cite{BelaDrin82,Semenov83}.

There are different constructions of the free commutative RB-algebra \cite{Rota68,Cartier72,GuoKeigher}.
In 2008, K. Ebrahimi-Fard and L. Guo constructed the free associative RB-algebra \cite{FardGuo07}.
In 2010, L.A. Bokut, Yu. Chen and X. Deng  \cite{BokutChen} found a~linear basis of the free associative RB-algebra
by means of Gr\"{o}bner---Shirshov technique. 
Linear bases of the free Lie RB-algebra were found in \cite{Gub2016,GubKol2016,Chen16}.

Pre-Lie algebras were introduced independently by E.B.~Vinberg \cite{Vinberg60} in 1960
and M.~Gerstenhaber \cite{Gerst63} in 1963. 
Pre-Lie algebras also known as left-symmetric algebras satisfy the identity 
$(x_1 x_2)x_3 - x_1(x_2 x_3) = (x_2 x_1) x_3 - x_2 (x_1 x_3)$.

In 2001, J.-L. Loday \cite{Dialg99} defined the dendriform dialgebra
(preassociative algebra) as a vector space endowed 
with two bilinear operations $\succ,\prec$ satisfying 
$$
\begin{gathered}
(x_1\succ x_2+x_1\prec x_2)\succ x_3 = x_1\succ (x_2 \succ x_3), \
(x_1\succ x_2)\prec x_3=x_1\succ(x_2\prec x_3), \\
x_1\prec(x_2\succ x_3+x_2\prec x_3)=(x_1\prec x_2)\prec x_3.
\end{gathered}
$$

In 1995, \cite{Loday95} J.-L. Loday also defined Zinbiel algebra
(precommutative algebra), on which the identity
$(x_1\succ x_2 + x_2\succ x_1)\succ x_3 = x_1\succ (x_2\succ x_3)$ holds. 
Any preassociative algebra
with the identity $x\succ y = y\prec x$ is a precommutative algebra
and under the product $x\cdot y = x\succ y - y\prec x$ is a pre-Lie algebra.

In 2004, dendriform trialgebra (postassociative algebra) was introduced~\cite{Trialg01}, 
i.e., an algebra with bilinear operations $\prec,\succ,\cdot$ satisfying seven certain axioms.
A space~$A$ with two bilinear operations $[,]$ and $\cdot$ is called a post-Lie algebra 
(Vallette, 2007 \cite{Vallette2007}) if $[,]$ is a Lie bracket and the next identities hold
$$
(x \cdot y) \cdot z - x \cdot (y \cdot z)
- (y \cdot x) \cdot z + y \cdot (x \cdot z) = [y,x]\cdot z, \quad
x\cdot [y,z] = [x \cdot y,z] + [y,x\cdot z].
$$

Let us explain the choice of terminology. 
Given a binary quadratic operad~$\mathcal P$, 
the defining identities for pre- and post-$\mathcal P$-algebras were found in \cite{BBGN2012}.
One can define the operad of pre- and post-$\mathcal P$-algebras
as $\mathcal P \bullet \mathrm{Pre}\Lie$ and
$\mathcal P \bullet \mathrm{Post}\Lie$ respectively.
Here $\mathrm{Pre}\Lie$ and $\mathrm{Post}\Lie$ denote 
the operads of pre-Lie algebras and post-Lie algebras respectively,
$\bullet$ denotes the black Manin product of operads \cite{GinzKapr}.
By pre- or postalgebra we will mean pre- or post-$\mathcal P$-algebra 
for some operad~$\mathcal P$.

In 2000, M. Aguiar \cite{Aguiar00} stated that any associative 
(commutative) algebra with a given Rota---Baxter operator~$R$
of weight zero under the operations $a \succ b = R(a)b$, $a\prec b = aR(b)$ 
is a preassociative (precommutative) algebra. In 2002, K. Ebrahimi-Fard~\cite{Fard02}
showed that any associative RB-algebra of nonzero weight $\lambda$ 
under the same two products $\succ$, $\prec$ and the third operation $a\cdot b = \lambda ab$ 
is a~postassociative algebra.
Analogous statement for postcommutative algebras was obtained by J.-L. Loday in 2007~\cite{Loday2007}.
The analogue of the Aguiar construction for the pair of pre-Lie algebras
and Lie RB-algebras of weight zero was stated in 2000 by M. Aguiar \cite{Aguiar00} 
and by I.Z. Golubchik, V.V. Sokolov~\cite{GolubchikSokolov}.
In 2010 \cite{BaiGuoNi10}, this construction for the pair of post-Lie algebras 
and Lie RB-algebras of nonzero weight was extended.

In 2013 \cite{BBGN2012}, the construction of M. Aguiar and K. Ebrahimi-Fard
was generalized for the case of arbitrary variety.

In 2008, the notion of universal enveloping RB-algebras
of pre- and postassociative algebras was introduced~\cite{FardGuo07}.
In~\cite{FardGuo07}, it was also proved that the universal enveloping of 
a free preassociative algebra is free.

In 2010, with the help of Gr\"{o}bner---Shirshov bases \cite{BokutChen},
Yu. Chen and Q.~Mo proved that any preassociative algebra
over the field of characteristic zero injectively embeds into
its universal enveloping RB-algebra \cite{Chen11}.

In 2013 \cite{GubKol2013}, given a variety $\Var$, it was proved that every pre-$\Var$-algebra
(post-$\Var$-algebra) injectively embeds into its universal enveloping $\Var$-RB-alge\-bra of weight 
$\lambda = 0$ ($\lambda\neq0$). Based on the last result, we formulate

{\bf Problem 1}.
Construct the universal enveloping RB-algebra of a pre- or post\-algebra.

In~\cite{Guo2011}, L.~Guo actually stated the following problem for $\Var = \As$:

{\bf Problem 2}.
Clarify if the pairs of varieties $(\RB\Var,\mathrm{pre}\Var)$
and $(\RB_\lambda\Var$, $\mathrm{post}\Var)$ for $\lambda\neq0$
are Poincar\'{e}---Birkhoff---Witt (PBW-) pairs in the sense of~\cite{PBW}.

Here $\RB\Var$ ($\RB_\lambda\Var$) denotes the variety of $\Var$-algebras
endowed with an RB-operator of (non)zero weight $\lambda$.

{\bf Problem 3}. (Guo et al., \cite{Guo2013})
Given a variety $\Var$ of algebras, whether the variety of RB-$\Var$-algebras is Schreier,
i.e., whether every subalgebra of the free algebra is free itself?

{\bf Problem 4}.
a) Prove that any pre-Lie (post-Lie) algebra injectively embeds into
its universal enveloping preassociative (postassociative) algebra.

b) Construct the universal enveloping preassociative 
(postassociative) algebra for given pre-Lie (post-Lie) algebra. 

The answer on Problem~4b will cover Problem~4a.
For pre-Lie algebras, Problem~4b and Problem~2 were stated in \cite{Kol2017}. 
The discussion of Problem~4 in the case of restricted 
pre-Lie algebras can be found in \cite{Dokas13}. 
The analogues of Problem~4 for Koszul-dual objects, di- and trialgebras, 
were solved in \cite{LP93,GubKol2014}.

Problems 1--3 were solved by author in commutative \cite{Gub2017Com}, 
associative \cite{Gub2017As}, and Lie \cite{Gub2017Lie} cases.
So, the question of L.~Guo \cite{Guo2011} is completely answered.

The main result of the current work is the affirmative answer on Problem~4a in postalgebra case.
In April 2018, the proof that the pair of varieties of pre-Lie
and preassociative algebras is a PBW-pair was announced~\cite{Dots}.

\section{Postalgebras}

{\em A postassociative algebra} is a linear space with three bilinear operations
$\cdot$, $\succ$, $\prec$ satisfying seven identities:
\begin{equation}\label{id:PostAs}
\begin{gathered}
(x \prec y) \prec z = x \prec (y \succ z + y \prec z + y\cdot z), \quad
(x \succ y) \prec z = x \succ (y \prec z), \\
(x \succ y + y\succ x + x\cdot y) \succ z = x \succ (y \succ z), \\
x\succ (y\cdot z) = (x\succ y) \cdot z,\quad
(x \prec y) \cdot z = x \cdot (y \succ z), \\
(x \cdot y) \prec z = x \cdot (y \prec z), \quad
(x \cdot y) \cdot z = x \cdot (y \cdot z).
\end{gathered}
\end{equation}

{\em A post-Lie algebra} is a vector space endowed with two bilinear products $[,]$ and $\cdot$,
the bracket $[,]$ is Lie, and the following identities are fulfilled:
$$
(x \cdot y) \cdot z - x \cdot (y \cdot z)
- (y \cdot x) \cdot z + y \cdot (x \cdot z) = [y,x]\cdot z, \quad
x\cdot [y,z] = [x \cdot y,z] + [y,x\cdot z].
$$

\section{Embedding of pre- and postalgebras into RB-algebras}

{\bf Theorem 1}~\cite{Aguiar00,BBGN2012,BaiGuoNi10,Fard02,GolubchikSokolov,Loday2007}.
Let $A$ be an RB-algebra of a variety $\Var$ and weight $\lambda=0$ ($\lambda\neq0$).
With respect to the operations
\begin{equation}\label{LodayToRB}
x\succ y = R(x)y,\quad x\prec y = xR(y)\quad (x\cdot y = \lambda xy)
\end{equation}
$A$ is a pre-$\Var$-algebra (post-$\Var$-algebra).

Denote the pre- and post-$\Var$-algebra obtained in Theorem~1 as $A^{(R)}_{\lambda}$. 

Given a pre-$\Var$-algebra $\langle C,\succ,\prec\rangle$, universal enveloping 
RB-$\Var$-algebra $U$ of $C$ is the universal algebra in the class of all RB-$\Var$-algebras
of weight zero such that there exists homomorphism from $C$ to $U^{(R)}_0$.
Analogously universal enveloping RB-$\Var$-algebra of a post-$\Var$-algebra is defined. 

{\bf Theorem 2}~\cite{GubKol2013}.
Any pre-$\Var$-algebra (post-$\Var$-algebra) could be embedded 
into its universal enveloping RB-algebra of the variety $\Var$ 
and weight $\lambda = 0$ ($\lambda\neq0$).

Let us consider the idea of the proof of Theorem~2 in the postalgebra case.
Suppose $\langle A,\succ,\prec,\cdot\rangle$ is a post-$\Var$-algebra.
Then the direct sum of two isomorphic copies of $A$,
the space $\hat A = A\oplus A'$, endowed with a binary operation
\begin{equation}\label{hat-construction}
a*b = a\succ b + a\prec b + a\cdot b, \
a*b' = (a\succ b)', \
a'*b = (a\prec b)', \
a'*b' = (\lambda a\cdot b)'
\end{equation}
for $a,b\in A$, is proved to be an algebra of the variety $\Var$.
Moreover, the map $R(a') = \lambda a$, $R(a) = -\lambda a$ is 
an RB-operator of weight~$\lambda$ on $\hat A$. 
The injective embedding of $A$ into $\hat A$ is given by $a\mapsto a'$, $a\in A$.
However, $\hat A$ is not a universal enveloping RB-algebra of~$A$.

\section{Embedding of post-Lie algebras into postassociative algebras}

Let $R\As\langle X\rangle$ denote the free associative algebra 
generated by a set $X$ with a linear map $R$ in the signature.
One can construct a linear basis of $R\As\langle X\rangle$ (see, e.g., \cite{Guo2013})
by induction. First, all elements from $S(X)$ lie in the basis.
Next, if we have basic elements $a_1,a_2,\ldots,a_k$, $k\geq1$, 
then the word $w_1R(a_1)w_2\ldots w_kR(a_k)w_{k+1}$ lies in the 
basis of $R\As\langle X\rangle$.
Here $w_1,\ldots,w_{k+1}\in S(X)\cup \emptyset$.
Let us denote the basis obtained as $RS(X)$.
Given a word~$u$ from $RS(X)$, the number of appeareances of the symbol~$R$
in~$u$ is denoted by~$\deg_R(u)$, $R$-degree of~$u$. 
We call any element from $RS(X)$ of the form $R(w)$ as $R$-letter.
By~$X_\infty$ we denote the union of $X$ and the set of all $R$-letters.
Given $u\in RS(X)$, define $\deg u$ (degree of $u$) as the length of $u$
in the alphabet~$X_\infty$.

Suppose that $X$ is a well-ordered set with respect to $<$.
Let us introduce by induction the deg-lex order on $S(X)$. 
Firstly, we compare two words $u$ and $v$ by the length:
$u < v$ if $|u|<|v|$. Secondly, when $|u| = |v|$, 
$u = x_i u'$, $v = x_j v'$, $x_i, x_j\in X$,
we have $u < v$ if either $x_i < x_j$ or $x_i = x_j$, $u'< v'$.
We compare two words $u$ and $v$ from $RS(X)$ 
by $R$-degree: $u<v$ if $\deg_R(u)<\deg_R(v)$. 
If $\deg_R(u) = \deg_R(v)$, we compare $u$ and $v$ in deg-lex order as words 
in the alphabet $X_\infty$. Here we define each $x$ from $X$ to 
be less than any $R$-letter and $R(a)<R(b)$ if and only if $a<b$.
Finally, we define an order $<$ on $R\As\langle X\rangle$,
extending it from basic elements.

Let $*$ be a symbol not containing in $X$.
By a $*$-bracketed word on $X$, we mean any basic word 
from $R\As\langle X\cup\{*\}\rangle$ with exactly one occurrence of $*$. 
The set of all $*$-bracketed words on $X$ is denoted by $RS^*(X)$.
For $q\in RS^*(X)$ and $u\in R\As\langle X\rangle$, we define
$q|_u$ as $ q|_{*\to u}$ to be the bracketed word obtained by 
replacing the letter $*$ in $q$ by $u$.

The order defined above is monomial, i.e., from $u < v$  
follows that $q|_u < q|_v$ for all $u,v\in RS(X)$ and $q\in RS^*(X)$.

Given $f\in R\As\langle X\rangle$, by $\bar{f}$ we mean the leading word in $f$.

{\bf Definition 1}~\cite{Guo2013}.
Let $f,g\in R\As\langle X\rangle$.
If there exist $\mu,\nu,w\in RS(X)$ such that 
$w = \bar{f} \mu = \nu \bar{g}$ with $\deg w<\deg(\bar{f})+\deg(\bar{g})$, 
then we define $(f,g)_w$ as $f\mu - \nu g$ and call it 
the {\it composition of intersection} of $f,g$ with respect to $(\mu,\nu)$.
If there exist $q\in RS^*(X)$ and $w\in RS(X)$ 
such that $w = \bar{f} = q|_{\bar{g}}$, then we define
$(f,g)^q_w$ as $f-q|_g$ and call it the {\it composition of inclusion} 
of $f,g$ with respect to $q$.

{\bf Definition 2}~\cite{Guo2013}.
Let $S$ be a subset of monic elements from $R\As\langle X\rangle$ 
and $w\in RS(X)$. 

(1) For $u,v\in R\As\langle X\rangle$, we call $u$ and $v$ 
congruent modulo $(S, w)$ and denote this by $u \equiv v \mod (S, w)$ 
if $u - v = \sum c_i q_i|_{s_i}$ with $c_i \in \Bbbk$, $q_i\in RS^*(X)$, 
$s_i\in S$ and $q_i|_{\overline{s_i}} < w$.

(2) For $f,g\in R\As\langle X\rangle$ and suitable 
$w,\mu,\nu$ or $q$ that give a composition of intersection
$(f,g)_w$ or a composition of inclusion $(f,g)^q_w$, 
the composition is called trivial modulo $(S, w)$ if
$(f,g)_w$ or $(f,g)^q_w \equiv 0 \mod (S, w)$.

(3) The set $S\subset R\As\langle X\rangle$ is 
called a {\it Gr\"{o}bner---Shirshov basis} if, for all $f,g\in S$, 
all compositions of intersection $(f,g)_w$ 
and all compositions of inclusion $(f,g)^q_w$ are trivial modulo $(S, w)$.

{\bf Theorem 3}~\cite{Guo2013}. 
Let $S$ be a set of monic elements in $R\As\langle X\rangle$, 
let $<$ be a monomial ordering on $RS(X)$ and let $Id(S)$ 
be the $R$-ideal of $R\As\langle X\rangle$ generated by $S$. 
If $S$ is a Gr\"{o}bner---Shirshov basis in $R\As\langle X\rangle$,
then $R\As\langle X\rangle = \Bbbk Irr(S)\oplus Id(S)$ 
where $Irr(S) = RS(X)\setminus \{q|_{\bar{s}}\mid q\in RS^*(X),s\in S\}$
and $Irr(S)$ is a linear basis of $R\As\langle X\rangle/Id(S)$.

Let $A$ be an associative algebra with an RB-operator $R$.
Then the algebra $A^{(-)}$ is a Lie RB-algebra
under the product $[x,y] = xy-yx$ and the same action of $R$. 
So, we state the analogue of Problem~1.4. 
How does the left adjoint functor from the category 
of Lie RB-algebras to the category of associative RB-algebras
(of the same weight) look like? Is it the embedding always injective?

Let $\hat{L} = L\oplus L'$ be exactly the Lie algebra 
with the RB-operator $R$ of weight~$-1$ constructed 
in the sketch of the proof of Theorem~2. 
Then $R(y_\alpha) = y_\alpha$, $R(x_\alpha) = 0$, and 
$X_\Lambda\cup Y_\Lambda = \{x_\alpha,\alpha\in\Lambda\}\cup \{y_\alpha,\alpha\in\Lambda\}$
is a linear basis of $\hat{L}$. Note that $\Span\{X\}$ and $\Span\{Y\}$
are Lie subalgebras of $\hat L$.
Suppose that $\Lambda$ is well-ordered set. 
Extend the order to $X_\Lambda\cup Y_\Lambda$ as follows: 
$t_\alpha<t_\beta$, $t\in\{x,y\}$, if and only if $\alpha<\beta$
and $y_\alpha<x_\beta$ for all $\alpha,\beta\in\Lambda$.

Consider the set $S$ of the following elements 
in $R\As\langle X_\Lambda\cup Y_\Lambda\rangle$:
\begin{equation}
\begin{gathered}\label{UnivRel}
x_{\alpha}x_\beta - x_\beta x_{\alpha} - [x_\alpha,x_\beta],\ 
y_{\alpha}y_\beta - y_\beta y_{\alpha} - [y_\alpha,y_\beta],\ \beta<\alpha,\\
x_{\alpha}y_\beta - y_\beta x_{\alpha} - [x_\alpha,y_\beta],
\end{gathered}
\end{equation}

\vspace{-0.7cm}

\begin{gather}
R(\vec{y}_\alpha) - \vec{y}_\alpha, \label{YRel} \\
R(\vec{x}_\alpha), \label{XRel} \\
\!\!\!R(R(a_1)\vec{y}_{\alpha_1}R(a_2){\ldots}R(a_{s})\vec{y}_{\alpha_s}R(a_{s+1}))
 {-} R(a_1)\vec{y}_{\alpha_1}R(a_2){\ldots} R(a_{s})\vec{y}_{\alpha_s}R(a_{s+1}), \label{LongRel} \\
R(R(a)\vec{x}_\alpha R(b)) - R(R(a)\vec{x}_\alpha b + a\vec{x}_\alpha R(b) - a\vec{x}_\alpha b), \label{XinMiddleRel} \\ 
R(R(a)\vec{x}_\alpha ) - R(a\vec{x}_\alpha ), \quad 
R(\vec{x}_\alpha R(b)) - R(\vec{x}_\alpha b), \label{XinRightRel} \\ 
R(a)R(b) - R(R(a)b + aR(b) - ab). \label{almostRB}
\end{gather}
Here $a,b,a_i$ are elements from $RS'(X_\Lambda\cup Y_\Lambda)
 = RS(X_\Lambda\cup Y_\Lambda)\setminus
 (S(X_\Lambda)\cup S(Y_\Lambda))$.
By $\vec{y}_{\alpha_i}$ or $\vec{x}_{\alpha}$ 
we mean a monomial from $\Bbbk[Y_\Lambda]$ or $\Bbbk[X_\Lambda]$ respectively.

{\bf Theorem 4}.
The set $S$ is a a Gr\"{o}bner---Shirshov basis in 
$R\As\langle X_\Lambda\cup Y_\Lambda\rangle$.

{\sc Proof}.
It is known that all compositions between two
elements from~\eqref{UnivRel} are trivial, as it is the method 
to construct the classical universal enveloping algebra.

Compositions of intersection between two expressions from \eqref{almostRB} 
are trivial by \cite{Guo2013}. 
Thus, all compositions of intersection which are not at the same time 
compositions of inclusion are trivial.

Let us consider a composition of inclusion between \eqref{XRel} and \eqref{almostRB}.
Let $\vec{x}_{\alpha}\in C(X)$, $a\in RS'(X_\Lambda\cup Y_\Lambda)$,
$w = R(a)R(\vec{x}_{\alpha})$. By the following, we get that the composition of inclusion is trivial:
\begin{gather*}
R(a)R(\vec{x}_{\alpha}) \mathop{\equiv}\limits^{\eqref{XRel}} 0  \mod (S,w), \\
R(a)R(\vec{x}_{\alpha}) \mathop{\equiv}\limits^{\eqref{almostRB}} 
 R(R(a)\vec{x}_{\alpha} + aR(\vec{x}_{\alpha}) - a\vec{x}_{\alpha}) \mathop{\equiv}\limits^{\eqref{XinRightRel}} 
 R(aR(\vec{x}_{\alpha})) \mathop{\equiv}\limits^{\eqref{XRel}} 0  \mod (S,w).
\end{gather*}

A composition of inclusion between \eqref{YRel} and \eqref{almostRB} is analogously trivial.

Compute a composition of inclusion between \eqref{LongRel} and \eqref{almostRB}.
For 
$\vec{y}_{\alpha_i}\in \Bbbk[Y]$,
$w = R(b)R(R(a_1)\vec{y}_{\alpha_1}\ldots \vec{y}_{\alpha_s}R(a_{s+1}))$,
$b,a_1,\ldots,a_{s+1}\in RS'(X_\Lambda\cup Y_\Lambda)$,
we have
\begin{multline*}\allowdisplaybreaks
R(b)R(R(a_1)\vec{y}_{\alpha_1}R(a_2)\ldots R(a_{s})\vec{y}_{\alpha_s}R(a_{s+1}))  \\
 \mathop{\equiv}\limits^{\eqref{almostRB}} 
R(R(b)R(a_1)\vec{y}_{\alpha_1}R(a_2)\ldots R(a_{s})\vec{y}_{\alpha_s}R(a_{s+1})) \\
 + R(bR(R(a_1)\vec{y}_{\alpha_1}R(a_2)\ldots R(a_{s})\vec{y}_{\alpha_s}R(a_{s+1}))) \\
 - R(bR(a_1)\vec{y}_{\alpha_1}R(a_2)\ldots R(a_{s})\vec{y}_{\alpha_s}R(a_{s+1})) \\
 \mathop{\equiv}\limits^{\eqref{LongRel}}
R(R(b)R(a_1)\vec{y}_{\alpha_1}R(a_2)\ldots R(a_{s})\vec{y}_{\alpha_s}R(a_{s+1})) \\
  \mathop{\equiv}\limits^{\eqref{almostRB}}
R(R(R(b)a_1+bR(a_1)-ba_1)\vec{y}_{\alpha_1}R(a_2)\ldots R(a_{s})\vec{y}_{\alpha_s}R(a_{s+1})) \\
  \mathop{\equiv}\limits^{\eqref{LongRel}}
R(R(b)a_1+bR(a_1)-ba_1)\vec{y}_{\alpha_1}R(a_2)\ldots R(a_{s})\vec{y}_{\alpha_s}R(a_{s+1}) \mod (S,w).   
\end{multline*}
\begin{multline*}
R(b)R(R(a_1)\vec{y}_{\alpha_1}R(a_2)\ldots R(a_{s})\vec{y}_{\alpha_s}R(a_{s+1})) \\
 \mathop{\equiv}\limits^{\eqref{LongRel}} 
R(b)R(a_1)\vec{y}_{\alpha_1}R(a_2)\ldots R(a_{s})\vec{y}_{\alpha_s}R(a_{s+1}) \\
 \mathop{\equiv}\limits^{\eqref{almostRB}} 
R(R(b)a_1+bR(a_1)-ba_1)\vec{y}_{\alpha_1}R(a_2)\ldots R(a_{s})\vec{y}_{\alpha_s}R(a_{s+1}) \mod (S,w). 
\end{multline*}
Thus, the corresponding composition of inclusion is trivial modulo $(S,w)$.

Let us calculate a composition of inclusion between \eqref{XinMiddleRel} and \eqref{almostRB}.
Let $\vec{x}_{\alpha}\in C(X)$, $a,b,c\in RS'(X_\Lambda\cup Y_\Lambda)$,
$w = R(a)R(R(b)\vec{x}_{\alpha}R(c))$. Modulo $(S,w)$ we have 
\begin{multline} \label{Comp:XinM-almRB-R}
R(a)R(R(b)\vec{x}_{\alpha}R(c)) \mathop{\equiv}\limits^{\eqref{XinMiddleRel}} 
 R(a)R(R(b)\vec{x}_{\alpha}c + b\vec{x}_{\alpha}R(c) - b\vec{x}_{\alpha}c) \\
 \mathop{\equiv}\limits^{\eqref{almostRB}} 
 R( R(a)R(b)\vec{x}_{\alpha}c + R(a)b\vec{x}_{\alpha}R(c) - R(a)b\vec{x}_{\alpha}c ) \\
 {+} R(aR(R(b)\vec{x}_{\alpha}c + b\vec{x}_{\alpha}R(c) - b\vec{x}_{\alpha}c)) 
 {-} R(aR(b)\vec{x}_{\alpha}c + ab\vec{x}_{\alpha}R(c) - ab\vec{x}_{\alpha}c ).
\end{multline}
For the first summand of the RHS of \eqref{Comp:XinM-almRB-R}, we obtain 
\begin{equation} \label{Comp:XinM-almRB-Ra}
R( R(a)R(b)\vec{x}_{\alpha}c) \mathop{\equiv}\limits^{\eqref{almostRB}} 
 R( R(R(a)b + aR(b)-ab)\vec{x}_{\alpha}c)  \mod (S,w).
\end{equation}

Also, modulo $(S,w)$, we have 
\begin{multline} \label{Comp:XinM-almRB-L} 
R(a)R(R(b)\vec{x}_{\alpha}R(c)) 
\mathop{\equiv}\limits^{\eqref{almostRB}} 
 R( R(a)R(b)\vec{x}_{\alpha}R(c) \\
 + aR(R(b)\vec{x}_{\alpha}R(c)) - aR(b)\vec{x}_{\alpha}R(c) )
\mathop{\equiv}\limits^{\eqref{XinMiddleRel}}  
 R( R(a)R(b)\vec{x}_{\alpha}R(c)) \\
 + R(aR(R(b)\vec{x}_{\alpha}c + b\vec{x}_{\alpha}R(c) - b\vec{x}_{\alpha}c)) - R(aR(b)\vec{x}_{\alpha}R(c)). 
\end{multline}
Further, the first summand of the RHS of \eqref{Comp:XinM-almRB-L}
is congruent modulo $(S,w)$ to
\begin{multline} \label{Comp:XinM-almRB-La}
R( R(a)R(b)\vec{x}_{\alpha}R(c)) \mathop{\equiv}\limits^{\eqref{almostRB}} 
R( R(R(a)b + aR(b)-ab)\vec{x}_{\alpha}R(c)) \\
\mathop{\equiv}\limits^{\eqref{XinMiddleRel}} 
 R( R(a)b\vec{x}_{\alpha}R(c) + aR(b)\vec{x}_{\alpha}R(c) - ab\vec{x}_{\alpha}R(c)) \\
 + R(R(R(a)b + aR(b)-ab)\vec{x}_{\alpha}c) 
 - R(R(a)b\vec{x}_{\alpha}c + aR(b)\vec{x}_{\alpha}c - ab\vec{x}_{\alpha}c).
\end{multline}

Note that substitution of \eqref{Comp:XinM-almRB-La}
into \eqref{Comp:XinM-almRB-L} gives exactly the same as 
the substitution of \eqref{Comp:XinM-almRB-Ra} into \eqref{Comp:XinM-almRB-R}.
Compositions of inclusion between \eqref{XinRightRel} and \eqref{almostRB} 
are trivial by similar calculations.
Theorem is proved.

{\bf Corollary~1}.
The quotient $A$ of $R\As\langle X_\Lambda\cup Y_\Lambda\rangle$ by $Id(S)$
is the universal enveloping associative RB-algebra of 
the Lie algebra $\hat{L}$ with the RB-operator $R$.
Moreover, $\hat{L}$ injectively embeds into $A^{(-)}$.

{\sc Proof}.
The RB-identity \eqref{RB} when at least one of arguments
lies in $S(X_\Lambda)\cup S(Y_\Lambda)$ follows
from \eqref{UnivRel}--\eqref{XRel}, \eqref{XinRightRel}
and the fact that the space $\Bbbk X_\Lambda$ is a subalgebra in $\hat{L}$.
Within \eqref{almostRB}, 
we get an asssociative RB-algebra structure on~$A$.
By \eqref{UnivRel}--\eqref{XRel}, 
we have that $A$ is enveloping of $\hat{L}$
for both: the Lie bracket $[,]$ and the action of $R$.

Let us prove that $A$ is the universal enveloping one.
At first, $A$ is generated by~$\hat{L}$. At second,
all elements from $S$ are identities in the universal enveloping
associative RB-algebra $U_{RB}(\hat{L})$. Indeed,
\eqref{UnivRel} are enveloping conditions for the product,
\eqref{almostRB} is the RB-identity.
Let us show by induction that~\eqref{YRel} and \eqref{XRel} 
follow from the enveloping conditions for the action of $R$ on $L'$
and RB-identity. Suppose that we have proved that 
$R(\vec{y}_\alpha) = \vec{y}_\alpha$ and
$R(\vec{y}_\beta) = \vec{y}_\beta$, then
\begin{equation}\label{YRel:proof}
\vec{y}_\alpha\vec{y}_\beta
 = R(\vec{y}_\alpha)R(\vec{y}_\beta)
 = R(R(\vec{y}_\alpha)\vec{y}_\beta 
  + \vec{y}_\alpha R(\vec{y}_\beta)-\vec{y}_\alpha\vec{y}_\beta)
 = R(\vec{y}_\alpha\vec{y}_\beta).
\end{equation}
Analogously we deduce $R(\vec{x}_\alpha) = 0$ for any $\vec{x}_\alpha\in S(X)$.

The relations~\eqref{XinRightRel} follow from~\eqref{XRel} and~\eqref{RB} immediately.
Now we state~\eqref{XinMiddleRel}:
\begin{multline*}
0 = R(a)R(\vec{x}_\alpha)R(b)
  = R(R(a)\vec{x}_\alpha-a\vec{x}_\alpha)R(b) \\
  = R( R(a)\vec{x}_\alpha R(b) -a\vec{x}_\alpha R(b) 
  + R(R(a)\vec{x}_\alpha-a\vec{x}_\alpha)b - R(a)\vec{x}_\alpha b + a\vec{x}_\alpha b  ) \\
  = R( R(a)\vec{x}_\alpha R(b) -a\vec{x}_\alpha R(b) - R(a)\vec{x}_\alpha b + a\vec{x}_\alpha b ).
\end{multline*}

Finally, let us prove~\eqref{LongRel} by induction on~$s$ analogously to~\eqref{YRel:proof}:
\begin{multline*}
R(a_1)\vec{y}_{\alpha_1}R(a_2)\ldots R(a_{s})\vec{y}_{\alpha_s}R(a_{s+1})
 = R(R(a_1)\vec{y}_{\alpha_1}\ldots R(a_{s}))R(\vec{y}_{\alpha_s}R(a_{s+1})) \\
 = R(R(a_1)\vec{y}_{\alpha_1}R(a_2){\ldots}R(a_{s})\vec{y}_{\alpha_s}R(a_{s+1})). 
\end{multline*}

Hence, the quotient of 
$R\As\langle X_\Lambda\cup Y_\Lambda\rangle$ by $Id(S)$ 
is the universal enveloping associative RB-algebra for $\hat{L}$.
By Theorems 3 and 4 we get the injectivity of embedding $\hat{L}$ into $A^{(-)}$.
Corollary is proved.

{\bf Remark}.
Analogously to Corollary~1, one can find the linear basis of the universal enveloping associative RB-algebra 
for a Lie algebra endowed with an RB-operator $R$ of nonzero weight $\lambda$ 
satisfying $R^2 = -\lambda R$.

{\bf Corollary~2}.
Any post-Lie algebra injectively embeds into its universal enveloping postassociative algebra.

{\sc Proof}.
Let $L$ be a post-Lie algebra. By Theorem~2, 
$L$ can be injectively embedded into $\hat{L}^{(R)}$
with the RB-operator $R$ of weight $-1$. 
Then, by Corollary~1, we embed the Lie RB-algebra $\hat{L}$ 
into the associative algebra $A$ with the RB-operator $P$. 
Thus, the subalgebra (in postalgebra sense) 
$T$ in $A^{(P)}$ generated by the set 
$L' = \Span\{x_\alpha - y_\alpha\mid \alpha\in \Lambda\}$ 
is an (injective) enveloping postassociative algebra 
of initial post-Lie algebra $L$.

The question whether the pair (postLie, postAs) is a PBW-pair is still open.

\section*{Acknowledgements}

The author expresses his gratitude to P. Kolesnikov for the useful comments and remarks.

This work was supported by the Austrian Science Foundation FWF grant P28079.

\noindent Vsevolod Gubarev \\
University of Vienna \\
Oskar-Morgenstern-Platz 1, 1090 Vienna, Austria \\
Sobolev Institute of mathematics \\
Acad. Koptyug ave. 4, 630090 Novosibirsk, Russia \\
e-mail: vsevolod.gubarev@univie.ac.at
\end{document}